\documentclass[leqno,11pt]{article}
\usepackage{latexsym}
\usepackage{amssymb}
\usepackage{amsfonts}
\usepackage{amsmath}
\usepackage{epsfig}
\usepackage{color}
\usepackage{pifont}
\usepackage{comment}
\usepackage{todonotes}

\allowdisplaybreaks
\makeatletter
\long\def\unmarkedfootnote#1{{\long\def\@makefntext##1{##1}\footnotetext{#1}}}
\makeatother

\newtheorem{definition}{Definition}[section]

\newtheorem{lemma}[definition]{Lemma}

\newtheorem{theorem}[definition]{Theorem}

\newtheorem{proposition}[definition]{Proposition}

\newtheorem{corollary}[definition]{Corollary}

\newtheorem{remark}[definition]{Remark}

\newtheorem{example}[definition]{Example}

\def\o{\Omega}

\def\m2{|\Omega | /2}
\def\M2{\frac{|\Omega |}{2}}
\def\u+{u_+^*}

\def\-p{\overline{p}}

\def\w0{{W_0^{1,p}(\Omega)}}
\def\i{{(0, |\Omega|)}}

\def\R{\mathbb R}
\def\N{\mathbb N}
\def\Z{\mathbb Z}
\def\ep{\varepsilon}
\def\rN{\mathbb R^N}
\def\bu{{\bf u}}
\newcommand\sn{{{\S}^{n-1}}}
\newcommand{\bm}[1]{{\boldsymbol{\rm #1}}}
\def\rn{{{\R}^n}}
\def\rN{{{\R}^N}}
\def\rNn{{{\R}^{Nn}}}

\newcommand{\hh}{{\cal H}^{n-1}}
\newcommand{\medint}{-\kern  -,395cm\int}
\newcommand{\medintinrigo}{-\kern  -,315cm\int}
\newcommand{\medelle}{-\kern  -,235cm L}
\newcommand{\medellenrigo}{-\kern  -,180cm L}
\newcommand{\qed}{\thinspace\null\nobreak\hfill
\hbox{\vbox{\kern-.2pt\hrule height.2pt
depth.2pt\kern-.2pt\kern-.2pt \hbox to1.8mm {\kern-.2pt\vrule
width.4pt \kern-.2pt\raise1.8mm\vbox to.2pt{} \lower0pt\vtop
to.2pt{}\hfil\kern-.2pt \vrule
width.4pt\kern-.2pt}\kern-.2pt\kern-.2pt \hrule height.2pt
depth.2pt \kern-.2pt}}\par\medbreak}

\title{Sobolev inequalities for the symmetric gradient  \\ in arbitrary domains
} 
\numberwithin{equation}{section}

\author{
  Andrea Cianchi\\
 {\it Dipartimento di Matematica e Informatica \lq\lq U. Dini", Universit\`a di Firenze}\\ {\it Viale Morgagni 67/A, 50134 Firenze, Italy} \\{\it  e-mail: cianchi@unifi.it}
\bigskip
\\
  Vladimir G. Maz'ya \\
  {\it   Department of Mathematics, Link\"oping University, SE-581
83 Link\"oping, Sweden}
  \\ and \\
{\it  RUDN University}\\ {\it
6 Miklukho-Maklay St, Moscow, 117198, Russia}
\\ {\it e-mail: vladimir.mazya@liu.se}
}

\date{}

\begin{document}
\maketitle

\begin{abstract}

A form of Sobolev inequalities for the symmetric gradient of vector-valued functions   is proposed, which allows for arbitrary ground domains  in $\rn$. In the relevant inequalities,   boundary regularity of domains is
replaced with information on boundary traces of trial functions. The  inequalities  so obtained exhibit the same exponents as in 
classical  inequalities for the full gradient of Sobolev functions, in regular domains. Furthermore, they involve  constants independent
of the geometry of the domain, and hence yield novel results
yet for smooth domains.  Our
approach relies upon a pointwise estimate for the functions in question via a Riesz potential of their symmetric gradient and an unconventional potential  depending on their boundary trace. 
\end{abstract}

\pagestyle{myheadings} \thispagestyle{plain}

 \unmarkedfootnote {
\par\noindent {\it Mathematics Subject
Classifications:} 46E35, 46E30.
\par\noindent {\it Key words and phrases: }  Sobolev inequalities, symmetric gradient, irregular domains, Riesz potential,
boundary traces, Lorentz spaces, Orlicz spaces.
}

\section{Introduction }\label{sec1}

Diverse mathematical models for physical phenomena, described by some vector-valued function $\bu: \Omega \to \rn$, depend on the derivatives of $\bu$ just through the symmetric part $\mathcal E \bu$ of its distributional gradient $\nabla \bu$. Here,    $\Omega$ is an open  set in $\rn$, $n \geq 2$,  and 
$$\mathcal E \bu = \tfrac 12\big(\nabla \bu + (\nabla \bu)^T\big)\,,$$
where  $(\nabla \bu)^T$ stands for the transpose matrix of $\nabla \bu$.
Instances in this connection are provided by the theory of non-Newtonian fluids, where $\bu$ represents the velocity of a fluid  \cite{AsMa, DL, fei3, FuS,
MalNRR, MR, MaRa}, and the theories of plasticity and nonlinear elasticity, where $\bu$ stands for the displacement of a body \cite{FuS, Kohn, Te}. The pertaining mathematical models amount to partial differential equations or variational integrals involving $\mathcal E \bu$. The regularity theory for this kind of problems, which is considerably less developed than the classical one, calls for the use of an ad hoc functional framework. 
\par This is provided by Sobolev type spaces, defined in analogy with the standard Sobolev spaces, where the role of $\nabla \bu$ is instead played by $\mathcal E \bu$. In particular, the homogeneous space $E^{1,p}(\Omega)$ is defined, for $p \geq 1$, as 
\begin{equation}\label{E}
E^{1,p}(\Omega) = \{\bu \in L^1_{\rm loc}(\Omega): \, |\mathcal E \bu| \in L^p(\Omega)\}.
\end{equation}
Its subspace $E^{1,p}_0(\Omega)$ is defined as the space of those 
functions that \lq\lq vanish" on $\partial \Omega$, in the sense that their  continuation   by $0$ outside $\Omega$ belongs to $E^{1,p}(\rn)$. 
%
Replacing $L^p(\Omega)$ in the definition of  $E^{1,p}(\Omega)$ by a more general Banach function space $X(\Omega)$  yields  a Sobolev type space for the symmetric gradient that will be denoted by $E^1X(\Omega)$. Namely,
\begin{equation}\label{EX}E^1X(\Omega) = \{\bu \in L^1_{\rm loc}(\Omega): \, |\mathcal E \bu| \in X(\Omega)\}.
\end{equation}
The subspace $E^1_0X(\Omega)$ can be defined accordingly. The need for these generalized spaces   arises, for example,  in   problems governed by nonlinearities of non-necessarily power type. They appear, for example, 
 in  the  Prandt-Eyring   fluids  \cite{BrDF,  E, FuS},
in  models  for  plastic  materials  with  logarithmic  hardening  \cite{FrS}, as well as in models for  the behavior of  fluids in certain liquid body armors \cite{HRJ, SMB, W} which are affected by exponential type nonlinearities.
\par Sobolev type inequalities for  spaces built upon the symmetric gradient are crucial in the regularity theory in question. 
 Sobolev inequalties for the space $E^{1,p}_0(\Omega)$ have exactly the same form as those for the standard Sobolev space $W^{1,p}_0(\Omega)$, in any   open set $\Omega $ of finite Lebesgue measure $\mathcal L^n (\Omega)$. In particular, when $p \in (1, \infty)$, this is a consequence of the Korn inequality, which tells us that, if $\bu \in E^{1,p}_0(\Omega)$, then $\bu$ is in fact weakly differentiable in $\Omega$, and $|\nabla \bu| \in L^p(\Omega)$. Moreover, 
\begin{equation}\label{korn}
\|\nabla \bu\|_{L^p(\Omega)} \leq c \|\mathcal E \bu\|_{L^p(\Omega)}
\end{equation}
for some constant $c=c(n,p, \mathcal L^n (\Omega))$ and every function $\bu \in E^{1,p}_0(\Omega)$ \cite{Fu1, Go1,Go2, Kohn, Korn, MM, Ne, Resh}. By contrast, this result fails in the endpoint cases when either    $p=1$ or $p=\infty$   \cite{CFM, LM, Or}. Still, Sobolev inequalities with the same target norms as for $W^{1,1}_0(\Omega)$ continue to hold in $E^{1,1}_0(\Omega)$ -- see \cite{strauss}, and the recent advances \cite{SV, Van}.
\\ If the zero boundary condition is dropped, namely if functions from $E^{1,p}(\Omega)$ are considered, and $\Omega$ is connected, a counterpart of inequality \eqref{korn} still holds, where the norm of $\nabla \bu$ on the left-hand side has to be replaced by the  $L^p$ distance of $\nabla \bu$ from the space of skew-symmetric $n\times n$ matrices. However, the resultant inequality requires suitable regularity assumptions on the domain $\Omega$. Lipschitz domains, or even John domains support this kind of Korn inequality \cite{DRS, DurMus04}.
In absence of a Korn type inequality, a reduction to the case of $W^{1,p}(\Omega)$ is not available, and a direct method for embeddings of $E^{1,p}(\Omega)$ is required. A technique based on potential estimates, reminiscent of the original proof by Sobolev, is developed in \cite{Camp}, and also applies  to less regular domains. 
\par In the present paper we offer a sort of  Sobolev type inequalities, involving the symmetric gradient, that requires no regularity on the domain at all. The point of view that will be adopted is that, in a Sobolev inequality for $\mathcal E \bu$, information on the integrability of the bounday trace of $\bu$  can serve as a substitute for boundary regularity of the ground domain $\Omega$, an assumption which is critical even in the usual situation when the full gradient $\nabla \bu$ is employed. 
\par To be more specific, the inequalities that will be dealt with have the form
\begin{equation}\label{basic}
\|\bu\|_{Y(\Omega, \mu)}  \leq c ( \|\mathcal E \bu\|_{X(\Omega )} +   \|\bu\|_{Z(\partial \Omega )}),
\end{equation}
where
$\|\cdot\|_{X(\Omega )}$ is a Banach function norm on $\Omega$ with
respect to  Lebesgue measure, $\|\cdot\|_{Y(\Omega ,
\mu)}$ is a Banach function norm with respect to a  possibly more
general  upper Alfhors regular measure $\mu$,  $\|\cdot\|_{Z(\partial \Omega )}$ is a Banach function norm on $\partial \Omega$ with respect to the $(n-1)$-dimensional Hausdorff measure $\hh$, and $c$ is a constant. 
Besides the arbitrariness of $\Omega$, a key feature of the inequalities to be established is  the fact  that the norms $\|\cdot\|_{X(\Omega )}$,  $\|\cdot\|_{Y(\Omega ,
\mu)}$ and  $\|\cdot\|_{Z(\partial \Omega )}$  in inequality \eqref{basic} are exactly the same as those appearing in a counterpart inequality on regular domains, with $\nabla \bu$ in  the place of $\mathcal E \bu$. Another distinctive trait is that the constant $c$ is independent of    the regularity of the domain $\Omega$. Thereby, under this respect, our conclusions are new even in the case of smooth domains.
\par As pointed out above, due to the lack of Korn type inequalities in irregular domains, inequalities like 
\eqref{basic}  cannot be derived via their counterparts for the full gradient
\begin{equation}\label{full}
\|\bu\|_{Y(\Omega, \mu)}  \leq c_1 \|\nabla  \bu\|_{X(\Omega )} +  c_2\|\bu\|_{Z(\partial \Omega )},
\end{equation}
that have recently been established, together
with higher-order versions, in \cite{cm_arbitrary}.
%
Another obstacle for such a derivation is that, even in regular domains $\Omega$, a Korn type inequality need not hold if $L^p(\Omega)$, with $p\in (1, \infty)$, is replaced by a more general Banach function norm $X(\Omega)$.  Let us notice that failure of the Korn inequality is not limited to the norms in $L^1(\Omega)$ and in $L^\infty(\Omega)$  \cite{BrD}.
 We refer to \cite{AcMi, BCD, BrD,  BMM, Ci2, DRS, Fu2} for positive results in the  case when the latter is a norm in an Orlicz space.
Let us incidentally mention that, in the special  situation when   $\mu = \mathcal L^n$,
$X =L^p$,  $Y = L^q$, and  $Z=L^r$, with
$1 \leq p < n$, $r\geq 1$ and $q = \min\{\tfrac{rn}{n-1},
\tfrac{np}{n-p}\}$, inequality \eqref{full} was established in
 \cite{Ma1960} via isoperimetric inequalities. The optimal constants $c_1$ and $c_2$ were also exhibited in that paper, for $p=1$, in the setting of scalar functions.  Mass transportation techniques have been exploited in  \cite{MaggiVillani1} to determine the optimal constants  when $1<p<n$. 
 Sharp constants in a parallel inequality, corresponding to 
 the borderline case when $p=n$, can be found in
\cite{MaggiVillani2}.
\par Our approach to inequality \eqref{basic} starts with   a fundamental pointwise estimate for functions $\bu \in E^{1,1}(\Omega)$, in terms of the Riesz potential  of order $1$ applied to $|\mathcal E \bu|$, plus a nonstandard potential   depending on the values of $\bu$ over $\partial \Omega$. This  estimate  is close in the spirit to the method of Sobolev, and   
its proof is inspired by arguments from \cite{Camp} and \cite{cm_arbitrary}.  It turns the problem of the validity of inequality \eqref{basic} into that of the boundedness of the apropos potential operators  between the function spaces $X(\Omega)$ and $Y(\Omega , \mu)$, and $Z(\partial \Omega)$ and $Y(\Omega , \mu)$, respectively. The new problem can be faced with the help of a reduction principle, which rests upon  an inequality, in rearrangement form, deduced from the fundamental pointwise inequality. This principle enables us to reduce the question of the boundedness of the relevant $n$-dimensional operators to the considerably simpler problem of the boundedness of one-dimensional Hardy type operators. With this tool at disposal, one can establish   inequalities of the form \eqref{basic} for various families of Banach funcion spaces $X$, $Y$ and $Z$. As examples, we present results for  Lebesgue, Lorentz and Zygmund spaces, and also obtain analogues of classical inequalities, such as the Yudovich-Pohozaev-Trudinger exponential inequality, as special instances in bordeline situations. A Rellich type compactness theorem is given as well.

\section{A pointwise estimate}\label{proofmain}

%

Our estimate for functions $\bu \in E^{1,1}(\Omega)$ at a point $x \in \Omega$ involves, loosely speaking, the value of the trace of $\bu$  at the first point on $\partial \Omega$   intercepted on each ray issued from  $x$. 
Since neither  regularity nor boundedness is a priori assumed on $\Omega$, traces of weakly differentiable functions on $\partial \Omega$  need not be defined. Our pointwise estimate an the other results of this paper will thus be stated for functions that are continuous in $\overline \Omega$. More precisely, we shall make use of the space
$$C_b(\overline \Omega) = \{\bu :  \hbox{$\bu$ is continuous in $\overline \Omega$ and has bounded support}\}.$$ 
Of course,  the norm inequalities that will be established for functions in  $E^1X(\Omega) \cap C_b(\overline \Omega)$, for some Banach function space $X(\Omega)$, continue to hold for  the closure of $C_b(\overline \Omega)$ in $E^1X(\Omega)$. Such a closure is known to agree with $E^1X(\Omega)$ itself if, for instance, $X=L^p$, and $\rn \setminus \overline \Omega$ satisfies the cone condition \cite[Proposition 1.3, Chapter 1]{Te}.
\par Let us begin our discussion with a few notations 
 and preliminary properties. Given any open set  $\Omega$ in $\rn$, with $n \geq 2$, and any point $x \in
 \Omega$,
 we set 
\begin{equation}\label{omegax}
\Omega _x = \{y \in \Omega: (1-t)x + ty \subset \Omega\,\,\,
\hbox{for every $t \in (0,1)$}\},
\end{equation}
and
\begin{equation}\label{partomegax}
(\partial \Omega )_x = \{y \in \partial \Omega: (1-t)x + ty \subset
\Omega\,\,\, \hbox{for every $t \in (0,1)$}\}.
\end{equation}
They are the largest subset of $\Omega$ and $\partial \Omega $,
respectively, which can be \lq\lq seen" from $x$.
 It is easily verified that $\Omega _x$ is
an open set. Furthermore, $(\partial
\Omega )_x$ is a Borel
%
set -- see  \cite[Proposition 3.1]{cm_arbitrary}.
\\
Next, define the sets
\begin{equation}\label{finite}
(\Omega \times \mathbb S ^{n-1})_0 = \{(x, \vartheta ) \in \Omega
\times \mathbb S ^{n-1}: x+t\vartheta \in \partial \Omega\,\,
\hbox{for some $t>0$}\},
\end{equation}
and
\begin{equation}\label{infinite}
(\Omega \times \mathbb S ^{n-1})_\infty = (\Omega \times \mathbb S
^{n-1}) \setminus (\Omega \times \mathbb S ^{n-1})_0\,,
\end{equation}
where $\mathbb  S^{n-1}$ denotes the $(n-1)$-dimensional unit sphere in $\rn$.
Clearly,
\begin{equation}\label{omegabound}
(\Omega \times \mathbb S ^{n-1})_0= \Omega \times \mathbb S ^{n-1}
\quad \hbox{if $\Omega $ is bounded.}
\end{equation}
 Let
 \begin{equation}\label{z}
 \zeta  : (\Omega \times \mathbb S ^{n-1})_0
\to \rn
\end{equation}
 be the function defined as
$$\zeta (x, \vartheta) = x + t \vartheta, \quad \hbox{where $t$ is
such that $x + t \vartheta \in (\partial \Omega) _x$}.$$
%
%
%
\\
 Given a function $g : \partial \Omega \to \R^m$,  $m \geq 1$, with compact
 support, we adopt the convention that  $g(\zeta (x, \vartheta))$ is defined for every
 $(x, \vartheta) \in \Omega \times \mathbb S ^{n-1}$, on extending it by $0$ on  $(\Omega \times \mathbb S ^{n-1})_\infty$; namely, we set
\begin{equation}\label{conv}
g(\zeta (x, \vartheta))= 0 \quad  \hbox{if $(x, \vartheta) \in
(\Omega \times \mathbb S ^{n-1})_\infty$.}
\end{equation}
%
Finaly, let us  introduce the function
\begin{equation}\label{abdef}
\mathfrak b: \Omega \times \mathbb S ^{n-1} \to (0,
\infty]
\end{equation}
given by
\begin{equation}\label{rbx}
\mathfrak  b(x, \vartheta ) = \begin{cases} |\zeta(x, \vartheta )-x|
& \hbox{if $(x, \vartheta ) \in  (\Omega \times \mathbb S
^{n-1})_0$,}
\\ \infty & \hbox{otherwise}.
\end{cases}
\end{equation}
%
One has that both $\zeta$ and  $\mathfrak b$ are Borel  functions \cite[Proposition 3.2]{cm_arbitrary}.
\par
We are now in a position to state our pointwise bound.

\begin{theorem}\label{intermest}
{\bf [Pointwise estimate]} Let $\Omega$ be any  open set in $\rn$,
$n \geq 2$. 
There exists a constant $C=C(n)$ such that
\begin{align}\label{hfund1ell}
|\bu  (x)| & \leq C
  \int _{\mathbb S^{n-1}}|\bu (\zeta (x, \vartheta ))|\,d\hh (\vartheta ) 
+ C \int _\o \frac{|\mathcal E \bu (y)|}{|x-y|^{n-1}}\, dy
\qquad \hbox{for $x \in \Omega$,}
\end{align}
%
for every function $\bu \in E^{1,1}(\Omega ) \cap C_b(\overline \Omega )$. Here,   convention
\eqref{conv} is adopted.
\end{theorem}

\begin{remark}\label{0}
{\rm Under the assumption that
\begin{equation}\label{zero}
 u=0 \quad \hbox{on $\partial
\Omega$, }
\end{equation}
inequality \eqref{hfund1ell} reduces to
\begin{align}\label{hfund1ell0}
|u(x)|  \leq C  \int _\o \frac{|\mathcal E \bu (y)|}{|x-y|^{n-1}}\, dy
\qquad \hbox{for a.e. $x \in \Omega$.}
\end{align}  }
\end{remark}

\noindent
{\bf Proof of Theorem \ref{intermest}}. 
Fix $x \in \Omega$. Given $\vartheta \in \mathbb S^{n-1}$,  consider the function 
$$[0, \mathfrak  b(x, \vartheta)] \ni t \mapsto \bu (x+t\vartheta )\cdot \vartheta\,.$$
We claim that, for $\hh$-a.e. $\vartheta \in \mathbb S^{n-1}$, this function   is locally absolutely continuous, and 
$$\frac{d}{dt}\bu (x+t\vartheta )\cdot \vartheta = \mathcal E \bu  (x+t\vartheta ) \vartheta \cdot \vartheta
%
 \quad \hbox{for a.e. $t \in [0, \mathfrak  b(x, \vartheta)]$.}$$
Note that this claim does not follow from standard properties of Sobolev functions, since $\bu$ need not belong to $W^{1,1}_{\rm loc}(\Omega)$, due to the failure of the Korn inequality for the $L^1$-norm. In order to prove it, let us begin by observing that, if $u$ is a smooth function, and we 
set $\bu = (u^1, \dots , u^n)$ and $\vartheta = (\vartheta _1, \dots , \vartheta _n)$, then
\begin{equation}\label{smooth} \bu (x+t\vartheta )\cdot \vartheta = \sum _{h=1}^n u^h (x+t\vartheta ) \vartheta _h\, \quad \hbox{and}
 \quad \mathcal E \bu  (x+t\vartheta ) \vartheta \cdot \vartheta = \sum _{h,k=1}^n u^h _{x_k}(x+t\vartheta ) \vartheta _h \vartheta _k\,.
\end{equation}
Indeed,
 $\sum _{h,k=1}^n u^h_{x_k} (x+t\vartheta )\vartheta _h \vartheta _k$ agrees with the quadratic form associated with the matrix $\nabla \bu (x+t\vartheta)$, evaluated at $\vartheta$. Therefore, it agrees with the quadratic form associated with the symmetric part $\mathcal E \bu (x+t\vartheta)$ of this matrix, evaluated at $\vartheta $, namely with $\mathcal E\bu (x+t\vartheta )\vartheta  \cdot \vartheta$.
\\
Next, let  $\Omega' $ be any smooth open set, starshaped with respect to $x$, and such that $\Omega ' \subset \subset \Omega$. 
Then the restriction of the function $\bu$ to $\Omega '$ can be extended to a function in $E^{1,1}\rn)$, still denoted by $\bu$,  with compact support -- see e.g. \cite[Remark 1.3, Chapter 2]{Te}. By \cite[Proposition 1.3, Chapter 1]{Te},  there exists a sequence of functions $\{\bu _m\} \subset C^\infty_0(\rn)$ such that $\bu _m \to \bu$ in $E^{1,1}(\rn)$, and $\bu_m \to \bu$ and $\mathcal E \bu_m \to \mathcal E \bu$ a.e. in $\rn$. In particular, $\{\bu _m\} $ is a Cauchy sequence in $E^{1,1}(\rn)$. Given $\varepsilon >0$, one can hence  make use of polar coordinates centered at $x$ and of \eqref{smooth} to deduce  that 
\begin{align}\label{sym20} 
\int _0^\infty & t^{n-1} \int _{\mathbb S^{n-1}}|\bu_m (x+t\vartheta )\cdot \vartheta - \bu_j (x+t\vartheta )\cdot \vartheta|\, d\hh (\vartheta)\, dt \\ \nonumber 
& \quad + \int _0^\infty t^{n-1} \int _{\mathbb S^{n-1}}\bigg|\frac{d}{dt}\bu_m (x+t\vartheta )\cdot \vartheta - \frac{d}{dt}\bu_j (x+t\vartheta )\cdot \vartheta\bigg|\, d\hh (\vartheta)\, dt
\\ \nonumber & = \int _0^\infty  t^{n-1} \int _{\mathbb S^{n-1}}|\bu_m (x+t\vartheta )\cdot \vartheta - \bu_j (x+t\vartheta )\cdot \vartheta|\, d\hh (\vartheta)\, dt \\ \nonumber
& \quad + \int _0^\infty t^{n-1} \int _{\mathbb S^{n-1}}\bigg|\sum _{h,k=1}^n (u_m)^h_{x_k} (x+t\vartheta )\vartheta _h \vartheta _k - \sum _{h,k=1}^n (u_j)^h_{x_k} (x+t\vartheta )\vartheta _h \vartheta _k\bigg|\, d\hh (\vartheta)\, dt
\\ \nonumber & = \int _0^\infty  t^{n-1} \int _{\mathbb S^{n-1}}|\bu_m (x+t\vartheta )\cdot \vartheta - \bu_j (x+t\vartheta )\cdot \vartheta|\, d\hh (\vartheta)\, dt \\ \nonumber
& \quad + \int _0^\infty t^{n-1} \int _{\mathbb S^{n-1}}\big|\big(\mathcal E \bu_m (x+t\vartheta) - \mathcal E \bu_j (x+t\vartheta )\big) \vartheta \cdot  \vartheta  \big|\, d\hh (\vartheta)\, dt
\\ \nonumber & \leq \int _0^\infty  t^{n-1} \int _{\mathbb S^{n-1}}|\bu_m (x+t\vartheta )\cdot \vartheta - \bu_j (x+t\vartheta )\cdot \vartheta|\, d\hh (\vartheta)\, dt \\ \nonumber
& \quad + C \int _0^\infty t^{n-1} \int _{\mathbb S^{n-1}}\big|\mathcal E\bu_m (x+t\vartheta) - \mathcal E \bu_j (x+t\vartheta ) \big|\, d\hh (\vartheta)\, dt
\\ \nonumber & \leq C \|\bu_m - \bu_j\|_{E^{1,1}(\rn)} < \varepsilon\,,
\end{align}
for some constant $C=C(n)$, 
provided that $m$ and $j$ are large enough. 
Hence, owing to Fubini's theorem,  for $\hh$-a.e. $\vartheta \in \mathbb S^{n-1}$,
\begin{align}\label{sym21} 
\int _0^\infty & t^{n-1}|\bu_m (x+t\vartheta )\cdot \vartheta - \bu_j (x+t\vartheta )\cdot \vartheta| \, dt \\ \nonumber 
& \quad + \int _0^\infty t^{n-1}  \bigg|\frac{d}{dt}\bu_m (x+t\vartheta )\cdot \vartheta - \frac{d}{dt}\bu_j (x+t\vartheta )\cdot \vartheta\bigg| \, dt < \varepsilon\,,
\end{align}
provided that $m$ and $j$ are large enough. Consequently, for a.e. $\vartheta \in \mathbb S^{n-1}$ the sequence  $\{\bu_m (x+ t \vartheta )\cdot \vartheta \}$ is a Cauchy sequence in $W^{1,1}(\delta , \infty)$ for any $\delta >0$. Hence, there exists a function   $\varphi _{x, \theta} \in W^{1,1}_{\rm loc}(0, \infty)$ such that 
\begin{equation}\label{sym22}
\bu_m (x+t\vartheta )\cdot \vartheta \to \varphi _{x, \theta}(t)
\end{equation}
in $W^{1,1}(\delta , \infty)$. Moreover, on taking a subsequence if necessary,
we may assume that
equation \eqref{sym22} holds for a.e. $t \in (0,\infty)$ and that
\begin{equation}\label{sym23}
\frac{d}{dt}\bu_m (x+t\vartheta )\cdot \vartheta \to \frac{d}{dt}\varphi _{x, \theta}(t)
\end{equation}
 for a.e. $t \in (0,\infty)$. On the other hand, since $\bu_m \to \bu$ and
$\mathcal E \bu_m \to \mathcal E \bu$ a.e. in $\rn$, by Fubini's theorem again, for $\hh$-a.e. $\vartheta \in \mathbb S^{n-1}$,
\begin{equation}\label{sym24}
\bu_m (x+t\vartheta )\cdot \vartheta \to  \bu (x+t\vartheta )\cdot \vartheta
\end{equation}
and 
\begin{equation}\label{sym25}
\frac{d}{dt}\bu_m (x+t\vartheta )\cdot \vartheta  = \mathcal E\bu_m (x+t\vartheta )\vartheta \cdot \vartheta \to \mathcal E\bu (x+t\vartheta )\vartheta \cdot \vartheta
\end{equation}
for a.e. $t \in (0, \infty)$. Equations \eqref{sym22}--\eqref{sym25} ensure that, for $\hh$-a.e. $\vartheta \in \mathbb S^{n-1}$, 
$$\varphi _{x, \theta}(t) = \bu (x+t\vartheta )\cdot \vartheta  \quad \hbox{and} \quad \frac{d}{dt}\varphi _{x, \theta}(t) = \mathcal E\bu (x+t\vartheta )\vartheta \cdot \vartheta$$ for a.e. $t \in (0, \infty)$. Hence, our claim follows, since $\Omega_x$ can be invaded by sets   $\Omega '$ as above.
\\ Now, let 
$\{\vartheta^1, \dots , \vartheta ^n\}\subset \mathbb S^{n-1}$ be a basis for $\rn$.  For fixed $i\in \{1, \dots, n\}$, define  the function $\varphi ^i : [0, \mathfrak  b(x, \vartheta^i )] \to \mathbb R$ as 
$$\varphi ^i (t) = \bu (x+t\vartheta ^i)\cdot \vartheta ^i \quad \hbox{for $t \in [0, \mathfrak  b(x, \vartheta^i )]$.}$$
In view of the property established above, for $\hh$-a.e. $\vartheta^1, \cdots , \vartheta ^n \in  \mathbb S^{n-1}$ forming a basis in  $\rn$, the function $\varphi^i$   is locally absolutely continuous, and 
$$\frac{d\varphi ^i (t)}{dt} = \mathcal E\bu (x+t\vartheta^i )\vartheta^i \cdot \vartheta^i
\quad \hbox{for $t \in [0, \mathfrak  b(x, \vartheta^i )]$.}$$
Thus, since $\varphi ^i (\mathfrak  b(x, \vartheta^i )) = \bu (\zeta (x, \vartheta^i))\cdot \vartheta ^i $,
\begin{align}\label{sym1}
 \bu (x)\cdot \vartheta ^i = \bu (\zeta (x, \vartheta^i))\cdot \vartheta ^i  - \int _0^{\mathfrak  b(x, \vartheta^i )} 
\mathcal E\bu (x+t\vartheta^i )\vartheta^i \cdot \vartheta^i
\, dt\,,
\end{align}
whence
\begin{align}\label{sym3}
| \bu (x)\cdot \vartheta ^i |\leq | \bu (\zeta (x, \vartheta^i))| + C \int _0^{\mathfrak  b(x, \vartheta^i )}  |\mathcal E \bu (x+t\vartheta ^i)| \, dt\,
\end{align}
for some constant $C=C(n)$.
We now exploit an argument from the proof of \cite[Lemma I.2]{Camp}. On setting 
$$a_{hk} = \sum _{i=1}^n \vartheta ^i_h \vartheta ^i_k \quad \quad \hbox{for $h, k =1, \dots ,n$,}$$
one has that
$$\sum _{i=1}^n(\bu (x)\cdot \vartheta ^i )^2 = \sum _{h,k=1}^n u^h(x) u^k(x) a_{hk}\,,$$
a quadratic form in $\bu (x)$ associated with the matrix ${\bf A} =\{a_{hk}\}$. Inasmuch as 
$\{\vartheta^1, \dots , \vartheta ^n\}$ is a basis in $\rn$, this quadratic form is positive definite. Furthermore, since the coefficients of its matrix ${\bf A}$ depend continuously on the unit vectors $\vartheta^1, \cdots , \vartheta ^n$, its smallest eigenvalue, which agrees with  the minimum of the quadratic form on the unit sphere, admits a positive lower bound $\lambda (K)$ as $(\vartheta^1, \dots , \vartheta ^n)$ ranges in a compact set $K \subset (\rn)^n$.
Therefore,
$$ \lambda (K)^{\frac 12}\, |\, \bu (x)|\leq  \Big(\sum _{i=1}^n(\bu (x)\cdot \vartheta ^i )^2\Big)^{\frac 12} \leq\sum _{i=1}^n|\bu (x)\cdot \vartheta ^i | \,.$$
As a consequence, inequality 	\eqref{sym3} implies that
\begin{align}\label{sym4}
\lambda (K)^{\frac 12} \, | \bu (x) |\leq \sum _{i=1}^n| \bu (\zeta (x, \vartheta^i))| + \sum _{i=1}^n\int _0^{\mathfrak  b(x, \vartheta^i )}  |\mathcal E \bu (x+t\vartheta ^i)| \, dt\,.
\end{align}
Next, set  $\eta = (\eta_1, \dots . \eta_{n-1})$,
$$Q = \{\eta\in \mathbb R^{n-1}: 0\leq \eta_i\leq 1,\,  i=1, \dots , n-1\},$$
and consider the functions $\phi_i : Q \to \mathbb R^n$, with $i=1, \dots , n$, defined as 
$$\phi _i (\eta)=\begin{cases} (\eta_1, 	\dots, \eta_{i-1}, \eta_i +1 , \eta_{i+1}, \dots , \eta_{n-1}, 0) &\quad \hbox{if $i=1, \dots , n-1$,}
\\  (\eta_1, 	\dots , \eta_{n-1}, 0) &\quad \hbox{if $ i=n$,}
\end{cases}
$$
for $\eta \in Q$.
Set $x_0= (1,1, \dots, 1) \in \rn$, and define the functions $\Phi _i : Q \to \mathbb S^{n-1}$, $ i=1, \dots , n$, as 
$$\Phi _i (\eta ) = \frac{\phi _i(\eta) -x_0}{|\phi _i(\eta) -x_0|} \quad \hbox{for $\eta \in Q$.}$$
For each index $i$, the function $\Phi _i : Q \to \Phi_i(Q)$ is a local coordinate system on $\mathbb S^{n-1}$. Moreover,
the set $\{\Phi_1(\eta), \dots , \Phi_n(\eta)\}$ is a basis of unit vectors in $\rn$ for each $\eta \in Q$, and the image of the map 
$$Q \ni \eta \mapsto (\Phi_1(\eta), \dots , \Phi_n(\eta)) \in (\rn)^n$$
is a compact subset of $(\rn)^n$. Thereby, an application of  inequality \eqref{sym4} with $\vartheta ^i= \Phi _i(\eta)$ implies that
\begin{align}\label{sym5}
c\, | \bu (x) |\leq \sum _{i=1}^n| \bu (\zeta (x, \Phi_i(\eta)))| + \sum _{i=1}^n\int _0^{\mathfrak  b(x, \Phi_i(\eta) )}  |\mathcal E \bu (x+t\Phi_i(\eta) )| \, dt\, \quad \hbox{for $\eta \in Q$,}
\end{align}
for a suitable constant $c=c(n)$.
A    change of variables ensures that
\begin{align}\label{sym6}
\int _Q| \bu (\zeta (x, \Phi_i(\eta)))|\, d\eta 
 \leq C \int _{\mathbb S^{n-1}} |\bu (\zeta (x, \vartheta))|\, d\mathcal H^{n-1}(\vartheta)
\end{align}
for some constant $C=C(n)$, and for $i=1, \dots , n$.
Analogously,
\begin{align}\label{sym7}
\int_Q
\int _0^{\mathfrak  b(x, \Phi_i(\eta) )} & |\mathcal E \bu (x+t\Phi_i(\eta) )| \, dt\, d\eta  \leq C
\int_{\mathbb S^{n-1}}
\int _0^{\mathfrak  b(x, \vartheta)}  |\mathcal E \bu (x+t\vartheta)| \, dt\, d\mathcal H^{n-1}(\vartheta)  \\ \nonumber & = C \int_{\mathbb S^{n-1}}
\int _0^{\mathfrak  b(x, \vartheta)}  t^{n-1}\frac{|\mathcal E \bu (x+t\vartheta)|}{t^{n-1}} \, dt\, d\mathcal H^{n-1}(\vartheta)  \leq C '\int_{\Omega_x}\frac{|\mathcal E \bu (y)|}{|x-y|^{n-1}}\, dy
\end{align}
for some constants $C=C(n)$, $C'=C'(n)$ and for $i=1, \dots , n$. Inequality \eqref{hfund1ell} follows on integrating inequality  \eqref{sym5} over $Q$, and exploiting equations \eqref{sym6} and \eqref{sym7}. \qed

\section{A rearrangement inequality and an ensuing reduction principle }\label{rearr}

The pointwise bound established in the previous section enables us
to derive an estimate, in rearrangement form, with respect to any $\alpha$-upper Ahlfors regular measure -- also called Frostman measure --
$\mu$ on $\Omega$,  with  $\alpha \in (n-1, n]$. Namely, a Borel measure $\mu$
such that
\begin{equation}\label{measure}
\mu (B_r(x) \cap \Omega ) \leq C_\mu r^\alpha\quad \hbox{for $x \in
\Omega$ and  $r>0$,}
\end{equation}
for some exponent  $\alpha \in (n-1, n]$  and some constant $C_\mu>0$. Here,
$B_r(x)$ denotes the ball, centered at $x$, with radius $r$.
\par
Given a  measure space $\mathcal R $, endowed with a
$\sigma$-finite, non-atomic measure $\nu$,  the decreasing rearrangement $\phi_\nu\sp*
: [0, \infty ) \to [0, \infty]$ of a $\nu$-measurable function $\phi
: \mathcal R \to \R$
%
%
%
 is defined as
$$
\phi_\nu\sp*(s)=\inf\{{t\geq 0}:\,\nu (\{|\phi|>t\})\leq s\}\quad
\textup{for}\ s\in [0,\infty).
$$
Although nonlinear, the operation of decreasing rearrangement has the property that
\begin{equation}\label{f+g}
(\phi + \psi)_\nu\sp*(s) \leq \phi _\nu\sp*(s/2) + \psi
_\nu\sp*(s/2)\quad \hbox{for $s \geq 0$,}
\end{equation}
for every measurable functions $\phi$ and $\psi$ on $\mathcal R$.
\\  Integrability properties of a function $\phi$ are preserved under the operation of decreasing rearrangement, since the functions $\phi$ and $\phi _\nu^*$ share the same distribution function. Indeed, 
$$\nu (\{|\phi|>t\}) = \mathcal L^1(\{\phi _\nu^* >t\}) \quad
\hbox{for every $t\geq 0$.}$$

%
%

\begin{theorem}\label{rearrintermest} {\bf [Rearrangement estimate]}
Let $\Omega$ be any  
open set in $\rn$, $n \geq 2$.  Assume that
$\mu$ is a Borel measure in $\Omega$ fulfilling \eqref{measure} for
some exponent $\alpha \in (n-1, n]$ and for some $C_\mu
>0$. Then there exists  constants $c=c(n)$ and $C=C(n, \alpha ,
C_\mu)$ such that
\begin{align}\label{rearresthell}
|\bu|_\mu^*(cs)  & \leq C\bigg[s^{-\frac
{n-1}{\alpha}}\int _0^{s^{\frac n\alpha}}|\mathcal E \bu|_{\mathcal
L^n}^*(r) \,dr + \int _{s^{\frac n\alpha}}^\infty r^{-\frac{n-1}n}
|\mathcal E \bu|_{\mathcal L^n}^*(r) \,dr
\\ \nonumber & \quad \quad \quad + s^{-\frac {n-1}{\alpha}}\int _0^{s^{\frac {n-1}\alpha}}
|\bu _{|\partial \Omega}|_{\hh}^*(r)\, dr\bigg]\quad
\quad \hbox{for $s>0$,}
\end{align}
for every function $\bu \in E^{1,1}(\Omega ) \cap C_b(\overline \Omega )$.
\end{theorem}

\smallskip
\noindent
{\bf Proof
 }. Let us denote by $I$ the classical Riesz
potential operator of order $1$, given by
\begin{equation}\label{I}
I f (x) = \int _\Omega \frac{f(y)}{|y-x|^{n-1}}\, dy
\quad \hbox{for   $x \in \Omega$,}
\end{equation}
at any $f \in L^1(\Omega)$. By \cite[Lemma 7.6]{cm_arbitrary}, 
 there exists a constant
$C=C(n, \alpha, C_\mu) $,    such that
\begin{align}\label{dic2}
(I f)_\mu^*(s) \leq C \bigg(s^{-\frac{n-1}{\alpha}} \int
_0^{s^{\frac{n}{\alpha}}} f^*_{\mathcal L^n}(r)\, dr + \int
_{s^{\frac{n}{\alpha}}}^\infty r^{-\frac {n-1}{n}}f
^*_{\mathcal L^n}(r)\, dr\bigg) \quad \hbox{for $s>0$,}
\end{align}
for every $f \in L^1(\Omega)$. Next,
define the operator $T$ as
\begin{equation}\label{T} Tg(x) = \int _{\mathbb S^{n-1}}|g(\zeta(x,
\vartheta ))|\,d\hh ( \vartheta ) \quad \hbox{ for $x \in \Omega$,}
\end{equation}
at any function Borel function $g:
\partial \Omega \to \R$. Here, and in what follows, we adopt convention
\eqref{conv}. Note that, owing to Fubini's theorem, $Tg$ is a
measurable function with respect to any Borel measure in $\Omega$. 
By \cite[Lemma 7.5]{cm_arbitrary},
there exists a
constant
 $C=C(n,\alpha , C_\mu)$ such that
\begin{align}\label{dic3}
(Tg)^*_\mu(s) \leq C s^{-\frac{n-1}{\alpha}}  \int
_0^{s^{\frac{n-1}{\alpha}}} g_{\hh}^*(r)\, dr \quad \hbox{for
$s>0$,}
\end{align}
for every Borel function $g : \partial \Omega \to \R$.
\\
With notations \eqref{I} and \eqref{T} in force, inequality \eqref{hfund1ell} takes the form
\begin{align}\label{dic1}
|\bu  (x)| & \leq C\big( T |\bu| (x) +  I  |\mathcal E \bu| (x)\big)
\qquad \hbox{for $x \in \Omega$.}
\end{align}
Hence, inequality \eqref{rearresthell}  follows via inequalities \eqref{dic2} and \eqref{dic3},
thanks to
property \eqref{f+g} of rearrangements. \qed

\bigskip
\par
Theorem \ref{rearrintermest} is the key step in a proof of the reduction principle contained in Theorem \ref{reduction} below. The latter enables one to derive norm inequalities of the form \eqref{basic} from corresponding one-dimensional inequalities for Hardy type operators dictated by the right-hand side of inequality \eqref{rearresthell}. Norms in rearrangement-invariant spaces are allowed in Theorem \ref{reduction}. Recalll that
%
%
%
%
a
rearrangement-invariant space $X(\mathcal R)$, on a measure space
$\mathcal R$ as above, is a Banach
function space (in the sense of Luxemburg) endowed with a norm
$\|\cdot \|_{X(\mathcal R)}$ such that
\begin{equation}\label{ri1}
\|\phi\|_{X(\mathcal R)} = \|\psi\|_{X(\mathcal R)}\quad
\hbox{whenever} \quad \phi^*_\nu = \psi^*_\nu.
\end{equation}
Every rearrangement-invariant space $X(\mathcal R)$ admits a
representation space $\overline X(0, \infty)$, namely another
rearrangement-invariant space  on $(0, \infty)$ such that
\begin{equation}\label{ri2}
\|\phi\|_{X(\mathcal R)} = \|\phi^*_\nu\|_{\overline X(0,
\infty)}\quad \hbox{for every $\phi \in X(\mathcal R)$.}
\end{equation}
In customary situations, an expression for the norm
$\|\cdot\|_{\overline X(0, \infty)}$ immediately follows from that
of $\|\cdot\|_{X(\mathcal R)}$. Lebesgue, Lorentz and Orlicz  spaces are  classical instances of rearrangement-invariant spaces. 
 We refer to \cite{BS} for a comprehensive account of
rearrangement-invariant spaces.

\begin{theorem}\label{reduction} {\bf [Reduction principle]}
Let $\Omega$ be any  open set in $\rn$, $n \geq 2$. Assume that
$\mu$ is a Borel measure in $\Omega$ fulfilling \eqref{measure} for some exponent
$\alpha \in (n-1, n]$ and for some constant $C_\mu$.  Let $X(\Omega)$, $Y(\Omega , \mu)$ and
 $Z(\partial \Omega)$ be rearrangement-invariant spaces such that
\begin{equation}\label{red1}
\left\| s^{-\frac {n-1}{\alpha}}\int _0^{s^{\frac
n\alpha}}\varphi(r) dr \right\|_{\overline Y(0, \infty )} \leq C
\|\varphi\|_{\overline X(0, \infty )},
\end{equation}
\begin{equation}\label{red2}
\left\| \int _{s^{\frac n\alpha}}^\infty r^{-\frac{n-1}n}
\varphi(r) dr \right\|_{\overline Y(0, \infty )} \leq C
\|\varphi\|_{\overline X(0, \infty )},
\end{equation}
\begin{equation}\label{red5}
\left\| s^{-\frac {n-1}{\alpha}}\int _0^{s^{\frac {n-1}\alpha}}
\varphi(r)\, dr \right\|_{\overline Y(0, \infty )} \leq C
\|\varphi\|_{\overline Z(0, \infty )},
\end{equation}
for some constant $C$, and for every non-increasing function
$\varphi : [0, \infty) \to [0, \infty)$.   Then there exists a constant $C'=C'(n, \alpha, C_\mu, C)$ such that
\begin{align}\label{red6}
\|\bu\|_{Y(\Omega, \mu)} & \leq C' \Big(\|\mathcal E \bu\|_{X(\Omega )} +  \|\bu\|_{Z(\partial \Omega )}\Big)
\end{align}
 for every  function $u \in E^1X(\Omega) \cap C_b(\overline \Omega)$.
\end{theorem}

\begin{remark}\label{remarkred}
{\rm Assumptions 
\eqref{red1}--\eqref{red5} of Theorem \ref{reduction}  can be weakened if either $\mu (\Omega )<
\infty$, or $\mathcal L^n (\Omega ) < \infty$, or  $\hh (\partial
\Omega ) < \infty$. Specifically: if $\mu (\Omega )< \infty$, it
suffices to assume that there exists $L \in (0, \infty)$ such that
inequalities \eqref{red1}--\eqref{red5} hold with the integral
operators multiplied by $\chi_{(0, L)}$ on the left-hand sides;  if
$\mathcal L^n (\Omega ) < \infty$, it suffices to assume that
inequalities \eqref{red1}--\eqref{red2} hold with $\varphi$ replaced
by $\varphi\chi_{(0, M)}$ for some $M \in (0, \infty)$; if $\hh
(\partial \Omega ) < \infty$, it suffices to assume that
inequality \eqref{red5} holds with $\varphi$ replaced
by $\varphi\chi_{(0, N)}$ for some $N \in (0, \infty)$. After these modifications in the assumptions,
inequality \eqref{red6} still holds, but with $C'$ depending also on
either on $L$ and $\mu (\Omega )$, or on $M$ and $\mathcal L^n
(\Omega )$, or on $N$ and $\hh (\partial \Omega ) < \infty$,
according to whether $\mu (\Omega )< \infty$, or $\mathcal L^n
(\Omega ) < \infty$, or $\hh (\partial \Omega ) < \infty$. }
\end{remark}

\medskip
\par\noindent
{\bf Proof of Theorem \ref{reduction}}. Let $u$ be any function as
in the statement. Inequalities \eqref{rearresthell} and
\eqref{red1}--\eqref{red5} imply that
\begin{align}\label{march1}
\big\| |\bu|_\mu^*(cs) \big\|_{\overline Y(0, \infty)} &
\leq C\Bigg[\bigg\|s^{-\frac {n-1}{\alpha}}\int _0^{s^{\frac
n\alpha}}|\mathcal E \bu|_{\mathcal L^n}^*(r) dr\bigg\|_{\overline
Y(0, \infty)} + \bigg\|\int _{s^{\frac n\alpha}}^\infty
r^{-\frac{n-1}n} |\mathcal E \bu|_{\mathcal L^n}^*(r)
dr\bigg\|_{\overline Y(0, \infty)}
\\ \nonumber & \quad  + \bigg\|s^{-\frac {n-1}{\alpha}}\int _0^{s^{\frac {n-1}\alpha}}
|\bu _{|\partial \Omega}|
_{\hh}^*(r)dr\bigg\|_{\overline Y(0, \infty)}\Bigg]
\\ \nonumber & \leq C' \Big[\big\||\mathcal E \bu|_{\mathcal L^n}^*(s)\|_ {\overline X(0, \infty)} + \big\||\bu _{|\partial \Omega}|
_{\hh}^*(s)\|_ {\overline Z(0, \infty)} \Big],
\end{align}
where $c=c(n)$, $C=C(n,\alpha , C_\mu)$ and $C'$ depends on $n,\alpha , C_\mu$ and on the constant $C$ appearing in \eqref{red1}--\eqref{red5}.
 On the other hand, a
property of rearrangement-invariant norms under dilations \cite[Proposition 5.11,
Chapter 3]{BS} tells us that the norm on the left-hand side of inequality
\eqref{march1} is bounded from below by $\min\big\{1, \tfrac 1c\big\} \big\|
|\bu|_\mu^*(s) \big\|_{\overline Y(0, \infty)}$. Thereby, one infers from
\eqref{march1} and \eqref{ri2}  that
\begin{align}\label{march2}
\|\bu\|_{Y(\Omega, \mu)} & \leq C  \big(\|\mathcal E \bu \|_{X(\Omega )} +   \|\bu \|_{Z(\partial \Omega )}\big)
\end{align}
for some constant $C$ depending on  $n,\alpha , C_\mu$ and on the constant $C$   in \eqref{red1}--\eqref{red5}. Inequality \eqref{red6}
follows.
 \qed

\section{Sobolev type inequalities}\label{sobineq}

In this section we exhibit a few  e inequalities  of the form \eqref{basic}, which  can be
established via Theorem \ref{reduction}. Inequalities involving  Lebesgue norms, Lorentz norms, and 
 Orlicz norms of exponential or power-logarithmic type, which
naturally come into play in borderline situations, will be presented. All inequalities are stated and proved for general
$\alpha$-upper Ahlfors regular measures 
$\mu$, with $\alpha \in (n-1, n]$. The  statements for the standard case of the Lebesgue measure can simply be obtained on setting $\alpha =n$. 
\par The target norm $\|\cdot\|_{Y(\Omega, \mu)}$ appearing in inequality \eqref{basic} depends on both the norm $\|\cdot\|_{X(\Omega)}$ and the norm $\|\cdot\|_{Z(\partial \Omega)}$. In the following discussion, given  $\|\cdot\|_{X(\Omega)}$, we  limit ourselves to considering the best possible norm $\|\cdot\|_{Z(\partial \Omega)}$, within a prescribed family of norms, in a boundary trace embedding for the space $E^1X(\Omega)$. The corresponding strongest possible norm $\|\cdot\|_{Y(\Omega, \mu)}$ in inequality \eqref{basic}  is  then exhibited. Of course, different choices of the norm $\|\cdot\|_{Z(\partial \Omega)}$ are possible for a given  $\|\cdot\|_{X(\Omega)}$. The optimal target norm  $\|\cdot\|_{Y(\Omega, \mu)}$ then depends on a balance between    $\|\cdot\|_{X(\Omega)}$ and $\|\cdot\|_{Z(\partial \Omega)}$.
\par We  premise some notations and  definitions 
 in connection with the norms appearing  in our results.
Let  $\mathcal R$ be a measure space equipped with a $\sigma$-finite, non-atomic measure $\nu$.
The Orlicz space
$L^A(\mathcal R)$ built upon a Young function $A: [0, \infty) \to
[0, \infty]$, namely  a left-continuous convex function which is
neither identically equal to $0$ nor to $\infty$, is a
rearrangement-invariant space equipped the Luxemburg norm given by
\begin{equation}\label{Orlicz}
\| \phi \|_{_{L^{A}(\mathcal R)}} =  \ \inf \Bigg\{ \lambda > 0 \ :
\ \int_{\mathcal R} A \bigg(\frac{|\phi (x)|}{\lambda} \bigg)\, dx \
\leq \ 1 \Bigg\}
\end{equation}
for  a measurable function $\phi$ in $\mathcal R$. 
The class of Orlicz spaces includes that of Lebesgue spaces, since
$L^{A}(\mathcal R)= L^{p}(\mathcal R)$ if $A(t)=t^p$ for $p \in [1,
\infty)$, and $L^{A}(\mathcal R)= L^{\infty}(\mathcal R)$ if
$A(t)=0$  for $t \in [0,1]$ and $A(t) = \infty$ for $t\in (1, \infty)$. Given $\sigma
>0$, we denote by $\exp L^\sigma (\mathcal R)$ and $\exp \exp L^\sigma (\mathcal R)$ the Orlicz spaces
built upon   Young functions equivalent to the functions $A(t)=e^{t^\sigma } -1$ and $A(t)= e^{e^{t^\sigma}} -e $, respectively, near inifnity. If $p>1$ and $\sigma \in \mathbb R$,  the notation
$L^p(\log L)^\sigma (\mathcal R)$ stands for the Orlicz space, also called Zygmund space, built upon a
Young function equivalent to the function $A(t) = t^p \log ^\sigma (1+t)$ near infinity. 
\\
The   Lorentz norms depend on two parameters, say $p$ and $q$. Assume that  either $1<p<\infty$
and  $1\leq q\leq\infty$, or $p=q=1$, or $p=q=\infty$. We define the functional
$\|\cdot\|_{L\sp{p,q}(\mathcal R)}$  by
$$
\|\phi\|_{L\sp{p,q}(\mathcal R)}=
\left\|s\sp{\frac{1}{p}-\frac{1}{q}}\phi^*_\nu(s)\right\|_{L\sp q(0,\infty)}
$$
for  a measurable function $\phi$ in $\mathcal R$. 
 Then $\|\cdot\|_{L\sp{p,q}(\mathcal R)}$ is equivalent to a~rearrangement-invariant  norm.
The corresponding space
$L\sp{p,q}(R)$ is called
Lorentz space.
\\
Suppose now that $1<q<\infty$ and $\nu (\mathcal R)< \infty$. The Lorentz-Zygmund 
space $L\sp{\infty,q;-1}(\mathcal R)$ is defined  via the rearrangement-invariant norm 
given by 
\begin{equation}\label{E:1.18}
\|\phi\|_{L\sp{\infty,q;-1}(\mathcal R)}=
\left\|s\sp{-\frac{1}{q}}\log \sp
{-1}\big(1+\tfrac{\nu (\mathcal R)}{s}\big) \phi^*_\nu(s)\right\|_{L\sp q(0,\nu (\mathcal R))}
\end{equation}
for    a measurable function $\phi$ in $\mathcal R$. 

\begin{theorem}\label{mainmeasure} {\bf [Subcritical Sobolev inequality]}
Let $\Omega$ be any  open set in $\rn$, $n \geq 2$. Assume that
$\mu$ is a measure in $\Omega$ fulfilling \eqref{measure} for some exponent
$\alpha \in (n-1, n]$ and for some constant $C_\mu$.   If $1 < p <
n$, then
 there
 exists a constant $C=C(n, p, \alpha, C_\mu)$ such that
\begin{align}\label{mainlebdisp}
\|\bu\|_{L^{\frac {\alpha p}{n-p}}(\Omega, \mu)} & \leq
C \Big(\|\mathcal E \bu\|_{L^p(\Omega )} +   
\|\bu\|_{L^{\frac{p(n-1)}{n-p}}(\partial \Omega )}\Big)
\end{align}
for every $u \in E^{1, p}(\Omega) \cap C_b(\overline \Omega)$.
\end{theorem}

\par\noindent
{\bf Proof}.  If $\mathcal R$ is any  $\sigma$-finite, non-atomic
 measure space, the space $L^p(0, \infty)$ is a representation space
of the Lebesgue space $L^p(\mathcal R)$. 
Inequality \eqref{mainlebdisp} then follows from Theorem \ref{reduction},
via standard one-dimensional Hardy type inequalities for Lebesgue
norms (see e.g. \cite[Section 1.3.2]{Mabook}). \qed

The next result tells us that, as in the classical Rellich theorem,
the Sobolev embedding corresponding to inequality
\eqref{mainlebdisp} is pre-compact  if the exponent $\tfrac {\alpha
p}{n- p}$ is replaced with any smaller one, and $\mu (\Omega) <
\infty$.

\begin{theorem}\label{maincompact} {\bf [Pre-compact Sobolev embedding]}
Let $\Omega$, $\mu$ and $p$ be as in Theorem
\ref{mainmeasure}. Assume, in addition, that $\mu (\Omega ) <
\infty$.  If $1 \leq q<\frac {\alpha p}{n-p}$, and $\{\bu_m\}$
is a bounded sequence in $E^{1,p}(\Omega) \cap L^{\frac{p(n-1)}{n-p}}(\partial \Omega)$,    then $\{\bu_m\}$ is a Cauchy sequence in $L^q(\Omega, \mu)$.
\end{theorem}

\medskip
\par\noindent
{\bf Proof}. Fix any $\varepsilon >0$, and 
choose a compact set $K\subset \Omega$ such that $\mu
(\Omega \setminus K)< \varepsilon$.  Let $\varrho \in C_0^\infty
(\Omega)$ be such that $0 \leq \varrho \leq 1$ and $\varrho = 1$ in
$K$. Thus, $K \subset  {\rm supp} (\varrho)$, the support of
$\varrho$, and hence
\begin{equation}\label{comp0}
\mu ({\rm supp} (1-\varrho )) \leq \mu(\Omega \setminus K)<
\varepsilon.
\end{equation}
 Let $\Omega '$ be an open set, with a smooth boundary,
satisfying ${\rm supp}(\varrho) \subset \Omega ' \subset \Omega$.
 Let $\{\bu_m\}$ be a bounded sequence in
$E^{1,p}(\Omega) \cap L^{\frac{p(n-1)}{n-p}}(\partial \Omega)$.
An application of Theorem \ref{mainmeasure},  with $\mu = \mathcal
L^n$, tells us that $\{\bu_m\}$ is also   bounded in $L^p(\Omega')$. As recalled in Section \ref{sec1}, thanks to a version of the Korn inequality, the space $E^{1,p}(\Omega') \cap L^p(\Omega ')$ agrees with the standard Sobolev space
$W^{1,p}(\Omega ')$, up to equivalent norms. By a weighted version of Rellich's
compactness theorem \cite[Theorem 1.4.6/1]{Mabook},
$\{\bu_m\}$ is a Cauchy sequence in $L^{q}(\Omega ', \mu)$,
and hence there exists $m_0 \in \N$ such that
\begin{equation}\label{comp1}
\|\bu_m - \bu_j\|_{L^q(\Omega ', \mu)}<\varepsilon
\end{equation}
if $m, j > m_0$. By H\"older's inequality,
\begin{align}\label{comp2}
\|(1-\varrho ) (\bu_m - \bu_j)\|_{L^q(\Omega, \mu)} &
\leq \|\bu_m - \bu_j\|_{L^{\frac {\alpha
p}{n-p}}(\Omega, \mu)} \mu ({\rm supp} (1-\varrho
))^{\frac{\alpha p - (n-p)q}{\alpha pq}}
\\ \nonumber
& \leq C \Big(\|\bu_m\|_{E^{1,p}(\Omega) \cap L^{\frac{p(n-1)}{n-p}}(\partial \Omega)} +\|\bu_j\|_{E^{1,p}(\Omega) \cap L^{\frac{p(n-1)}{n-p}}(\partial \Omega)} \Big)
\varepsilon^{\frac{\alpha p
- (n-p)q}{\alpha pq}}\\
\nonumber & \leq C'\varepsilon^{\frac{\alpha p - (n-p)q}{\alpha
pq}}
\end{align}
for some constants $C$ and $C'$ independent of $m$ and $j$. Inequalities
\eqref{comp1} and \eqref{comp2} tell us that
\begin{equation}\label{comp3}
\|\bu_m - \bu_j\|_{L^q(\Omega , \mu)} \leq \|\bu_m - \bu_j\|_{L^q(\Omega ', \mu)} +  \|(1-\varrho
)(\bu_m - \bu_j)\|_{L^q(\Omega, \mu)} \leq
\varepsilon + C'\varepsilon^{\frac{\alpha p - (n-p)q}{\alpha pq}}
\end{equation}
if $m, j > m_0$. Owing to the arbitrariness of $\varepsilon$,
inequality \eqref{comp3} implies that $\{\bu_m\}$ is a
Cauchy sequence in $L^{q}(\Omega , \mu)$. \qed

\bigskip
\par
 The next statement concernes the borderline exponent $p= n$, which is not included in
Theorem \ref{mainmeasure}.  It consists of two inequalities. The former is
  a version,  for the symmetric gradient, of the Yudovich-Pohozaev-Trudinger inequality,
in arbitrary domains with an upper Ahlfors regular measure, and involves an Orlicz norm of exponential type. The latter amonts to a slight improvement of the former, in that it allows for a stronger norm in a Lorentz-Zygmund space. It can be regarded as a counterpart in the  present framework of a result independently obtained by   Hansson and Brezis-Weinger in the classical setting. Interestingly, the norms coming into play in this limiting situation are independent of the exponent $\alpha$ in \eqref{measure}, and hence they are the same as for the Lebesgue measure. They also agree with those appearing in the parallel inequalities for the full gradient.  By contrast, the exponent $\alpha$ affects the constants in the inequalities in question.  Let us incidentally mention that the sharp constant in the exponential Yudovich-Pohozaev-Trudinger inequality, for the  standard gradient,  is detected in the paper \cite{cianchimoser}, which extends to Ahlfors regular measures a classical inequality by Moser.

\begin{theorem}\label{trudmeasure} {\bf [Critical Sobolev inequality]}
Let $\Omega$ and $\mu$ be as in Theorem \ref{mainmeasure} . Assume,
in addition, that $\mathcal L^n (\Omega)< \infty$, $\mu (\Omega)<
\infty$ and $ \hh (\partial \Omega) < \infty$. 
Then there
 exists a constant $C=C(n,\alpha, C_\mu, \mathcal L^n (\Omega), \mu(\Omega), \hh (\partial \Omega))$ such that
\begin{align}\label{truddisp}
\|\bu\|_{\exp L^{\frac {n}{n-1}}(\Omega, \mu)}  \leq
C\Big( \|\mathcal E \bu\|_{L^{n}(\Omega )}
  +
 \|\bu\|_{\exp L^{\frac
{n}{n-1}}(\partial \Omega)}\Big)
\end{align}
for every $u \in E^{1, n}(\Omega) \cap C_b(\overline \Omega )$.
\\ Moreover, there exists a constant as above such that
\begin{align}\label{BWp}
\|\bu\|_{L^{\infty; n,-1}(\Omega, \mu)}  \leq
C\Big( \|\mathcal E \bu\|_{L^{n}(\Omega )}
  +
 \|\bu\|_{L^{\infty; n,-1}(\partial \Omega)}\Big)
\end{align}
for every $u \in E^{1, n}(\Omega) \cap C_b(\overline \Omega )$.
\end{theorem}

\medskip
\par\noindent
{\bf Proof}. If $\mathcal R$ is a
finite measure space, then the norm of a function $\phi$ in the
Orlicz space $\exp L^\sigma (\mathcal R)$, with $\sigma
>0$, is equivalent, up to multiplicative constants depending on
$\sigma$ and $\nu (\mathcal R)$, to the functional
$$
\Big\|\big(1+ \log \tfrac {\nu(\mathcal R)}{s}\big)^{-\frac 1\sigma}
\phi _\nu ^*(s)\Big\|_{L^\infty (0, \nu (\mathcal R))},$$
see e.g. \cite[Lemma 6.12, Chapter 4]{BS}. Now, one can verify that the function 
\begin{equation}\label{mono1}(0, \infty) \ni s \mapsto s^{-\frac {n-1}{\alpha}}\int _0^{s^{\frac
n\alpha}}\varphi(r) dr   + \int _{s^{\frac n\alpha}}^\infty r^{-\frac{n-1}n}
\varphi(r) dr\,,
\end{equation}
is non-decreasing. Also, if $\varphi$ is non-increasing, then the function
\begin{equation}\label{mono2}
(0, \infty) \ni s \mapsto s^{-\frac {n-1}{\alpha}}\int _0^{s^{\frac
{n-1}\alpha}}\varphi(r) dr \,,
\end{equation}
is non-increasing as well.
By
Theorem \ref{reduction} and Remark \ref{remarkred}, the proof of inequality
 \eqref{truddisp} is thus reduced to showing that
\begin{equation}\label{prooftrud1}
\left\| s^{-\frac {n-1}{\alpha}}\big(1+ \log \tfrac {\mu
(\Omega)}s\big)^{-\frac{n-1}{n}}\int _0^{s^{\frac
n\alpha}}\varphi(r) dr \right\|_{L^{\infty}(0, \mu (\Omega) )} \leq
C \|\varphi\|_{L^{n}(0, \mathcal L^n (\Omega) )},
\end{equation}
\begin{equation}\label{prooftrud2}
\left\| \big(1+ \log \tfrac {\mu (\Omega)}s\big)^{-\frac{n-1}{n}}\int _{s^{\frac n\alpha}}^\infty r^{-\frac{n-1}n}
\varphi(r) dr \right\|_{L^{\infty}(0, \mu (\Omega) )} \leq C
\|\varphi\|_{L^{n }(0, \mathcal L^n(\Omega) )},
\end{equation}
for every non-increasing function $\varphi: (0, \infty) \to [0,
\infty)$ with support in $(0, \mathcal L^n(\Omega))$, and
\begin{multline}\label{prooftrud7}
\left\| s^{-\frac {n-1}{\alpha}}\big(1+ \log \tfrac {\mu
(\Omega)}s\big)^{-\frac{n-1}{n}}\int _0^{s^{\frac {n-1}\alpha}}
\varphi(r)\, dr \right\|_{L^{\infty}(0, \mu (\Omega) )} \\ \leq C
\Big\|\big(1+ \log \tfrac {\hh(\partial \Omega)}s\big)^{-\frac{n-1}{n}}\varphi(s)\Big\|_{L^{\infty}(0, \hh(\partial \Omega) )},
\end{multline}
for some constant $C$ and every non-increasing function for every non-increasing function $\varphi: (0, \infty) \to [0,
\infty)$ with support in $(0, \hh(\partial \Omega))$
\\
Inequalities \eqref{prooftrud1}--\eqref{prooftrud7} can be verified via classical characterizations of  Hardy type inequalities in weighted Lebesgue spaces \cite[Section
1.3.2]{Mabook}. 
\\ Let us now consider inequality \eqref{BWp}. By
Theorem \ref{reduction} and Remark \ref{remarkred} again, it suffices to show that
\begin{equation}\label{BW1}
\left\| s^{-\frac {n-1}{\alpha}}\int _0^{s^{\frac
n\alpha}}\varphi(r) dr   + \int _{s^{\frac n\alpha}}^\infty r^{-\frac{n-1}n}
\varphi(r) dr  \right\|_{L^{\infty; n,-1}(0, \mu (\Omega) )} \leq
C \|\varphi\|_{L^{n}(0, \mathcal L^n (\Omega) )},
\end{equation}
for some constant $C$ and 
for every non-increasing function $\varphi: (0, \infty) \to [0,
\infty)$ with support in $(0, \mathcal L^n(\Omega))$, and
\begin{equation}\label{BW3}
\left\| s^{-\frac {n-1}{\alpha}}\int _0^{s^{\frac {n-1}\alpha}}
\varphi(r)\, dr \right\|_{L^{\infty; n,-1}(0, \mu (\Omega) )} \leq C
\|\varphi\|_{L^{\infty; n,-1}(0, \hh(\partial \Omega) )},
\end{equation}
for some constant $C$ and every non-increasing function for every non-increasing function $\varphi: (0, \infty) \to [0,
\infty)$ with support in $(0, \hh(\partial \Omega))$. 
As a consequence of the monotonicy of the functions in \eqref{mono1} and \eqref{mono2}, inequality \eqref{BW1} is equivalent to the couple of inequalities 
\begin{equation}\label{BW4}
\left\| s^{-\frac {n-1}{\alpha}-\frac 1n}  \big(1+ \log \tfrac {\mu (\Omega)}s\big)^{-1}  \int _0^{s^{\frac
n\alpha}}\varphi(r) dr    \right\|_{L^{n}(0, \mu (\Omega) )} \leq
C \|\varphi\|_{L^{n}(0, \mathcal L^n (\Omega) )},
\end{equation}
and 
\begin{equation}\label{BW5}
\left\| \int _{s^{\frac n\alpha}}^\infty r^{-\frac{n-1}n}
\varphi(r) dr  \right\|_{L^{\infty; n,-1}(0, \mu (\Omega) )} \leq
C \|\varphi\|_{L^{n}(0, \mathcal L^n (\Omega) )}
\end{equation}
for some constant $C$ and 
for every non-increasing function $\varphi: (0, \infty) \to [0,
\infty)$ with support in $(0, \mathcal L^n(\Omega))$, 
and inequality \eqref{BW3} is equivalent to the inequality
%
\begin{equation}\label{BW6}
\left\| s^{-\frac {n-1}{\alpha}-\frac 1n} \big(1+ \log \tfrac {\mu (\Omega)}s\big)^{-1} \int _0^{s^{\frac {n-1}\alpha}}
\varphi(r)\, dr \right\|_{L^{n}(0, \mu (\Omega) )} \leq C
\|\varphi\|_{L^{\infty; n,-1}(0, \hh(\partial \Omega) )},
\end{equation}
for some constant $C$ and every non-increasing function for every non-increasing function $\varphi: (0, \infty) \to [0,
\infty)$ with support in $(0, \hh(\partial \Omega))$. Inequalities \eqref{BW4} and \eqref{BW6} can be derived as special cases of  Hardy type inequalities in weighted Lebesgue spaces  \cite[Section
1.3.2]{Mabook}. The proof of inequality \eqref{BW5} is subtler, and makes use of the fact that non-increasing trial functions $\varphi$ are considered. It follows via the same proof as  (one of the cases)   of \cite[Theorem 5.1]{CP_TAMS}. \qed

\bigskip
\par
The super-critical regime,  corresponding to the case when $p> n$, is the subject
of the following theorem.

\begin{theorem}\label{inf} {\bf [Super--critical Sobolev inequality]}
Let $\Omega$ be a  open set in $\rn$, $n \geq 2$, such that
 $\mathcal L^n (\Omega)< \infty$  and $ \hh (\partial \Omega) < \infty$.   If $p
>n$, then
 there
 exists a  constant $C=C(n, p,   \mathcal L^n (\Omega),  \hh (\partial
 \Omega))$
  such that
\begin{align}\label{infdisp}
\|\bu\|_{L^\infty (\Omega )} & \leq C \Big(\|\mathcal E \bu\|_{L^{p}(\Omega )}
 +     \|\bu\|_{L^\infty(\partial \Omega)}\Big)
\end{align}
for every $u \in E^{1,p}(\Omega) \cap C_b(\overline \Omega)$.
\end{theorem}

\medskip
\par\noindent
{\bf Proof of Theorem \ref{inf}}. Inequality \eqref{infdisp} follows
from Theorem \ref{reduction} and Remark \ref{remarkred}, via
weighted Hardy type inequalities  (\cite[Section 1.3.2]{Mabook}).
 \qed

\medskip
The next  result, contained in Theorem \ref{Lorentz}, concerns inequalities for functions whose  symmetric gradient belongs to a Lorentz space $L^{p,q}(\Omega)$. It extends   Theorems \ref{mainmeasure}, \ref{trudmeasure} and \ref{inf}, since $L^{p,p}(\Omega) = L^p(\Omega)$.  In fact, the conclusion of Theorem \ref{mainmeasure} 
is also augmented by the result of Part (i) of Theorem \ref{Lorentz}, with che choice   $q=p$. Actually, the Lorentz space $L^{\frac {\alpha p}{n-p},p}(\Omega, \mu)$, which is obtained as a target space by Theorem \ref{Lorentz}, is strictly contained in the Lebesgue space $L^{\frac {\alpha p}{n-p}}(\Omega, \mu)$ given by Theorem \ref{mainmeasure}.

\begin{theorem}\label{Lorentz} {\bf [Lorentz--Sobolev inequalities]}
Let $\Omega$ be any  open set in $\rn$, $n \geq 2$. Assume that
$\mu$ is a measure in $\Omega$ fulfilling \eqref{measure} for some exponent 
$\alpha \in (n-1, n]$  and for some constant $C_\mu$.  
\\ (i) Assume  that$1 < p <
n$ and $1\leq q \leq \infty$. Then
 there
 exists a constant $C=C(n, p, q, \alpha, C_\mu)$ such that
\begin{align}\label{mainlor}
\|\bu\|_{L^{\frac {\alpha p}{n-p},q}(\Omega, \mu)} & \leq
C \Big(\|\mathcal E \bu\|_{L^{p,q}(\Omega )} +   
\|\bu\|_{L^{\frac{p(n-1)}{n-p},q}(\partial \Omega )}\Big)
\end{align}
for every $u \in E^{1}L^{p,q}(\Omega) \cap C_b(\overline \Omega)$.
\\ (ii) Assume that $p=n$ and $q >1$, and that  $\mathcal L^n (\Omega)< \infty$, $\mu (\Omega)<
\infty$ and $ \hh (\partial \Omega) < \infty$. 
Then there
 exists a constant $C=C(n, q, \alpha, C_\mu, \mathcal L^n (\Omega), \mu(\Omega), \hh (\partial \Omega))$ such that
 \begin{align}\label{BWlor}
\|\bu\|_{L^{\infty; q,-1}(\Omega, \mu)}  \leq
C\Big( \|\mathcal E \bu\|_{L^{n,q}(\Omega )}
  +
 \|\bu\|_{L^{\infty; q,-1}(\Omega)}\Big)
\end{align}
for every $u \in E^{1}L^{n,q}(\Omega) \cap C_b(\overline \Omega )$.
\\ (iii) Assume that either $p=n$ and $q =1$, or $p>n$ and $1\leq q\leq \infty$, and that  $\mathcal L^n (\Omega)< \infty$, $\mu (\Omega)<
\infty$ and $ \hh (\partial \Omega) < \infty$. 
Then there
 exists a constant $C=C(n, p, q, \alpha, C_\mu, \mathcal L^n (\Omega), \mu(\Omega), \hh (\partial \Omega))$ such that
\begin{align}\label{inflor}
\|\bu\|_{L^\infty (\Omega )} & \leq C \Big(\|\mathcal E \bu\|_{L^{p,q}(\Omega )}
 +     \|\bu\|_{L^\infty(\partial \Omega)}\Big)
\end{align}
for every $u \in E^{1}L^{p,q}(\Omega) \cap C_b(\overline \Omega)$.
\end{theorem}

\begin{remark}\label{BW}
{\rm Under the assumptions of Theorem \ref{Lorentz}, Part (ii), an inequality (slightly weaker than \eqref{BWlor}) involving Orlicz norms of exponential type can be shown to hold. It  extends    \eqref{truddisp},  and tells us that there exists a constant $C=C(n, q, \alpha, C_\mu, \mathcal L^n (\Omega), \mu(\Omega), \hh (\partial \Omega))$, such that
\begin{align}\label{trudlor}
\|\bu\|_{\exp L^{\frac {q}{q-1}}(\Omega, \mu)}  \leq
C\Big( \|\mathcal E \bu\|_{L^{n,q}(\Omega )}
  +
 \|\bu\|_{\exp L^{\frac
{q}{q-1}}(\partial \Omega)}\Big)
\end{align}
for every $u \in E^{1}L^{n,q}(\Omega) \cap C_b(\overline \Omega )$.
}
\end{remark}

\medskip
\par\noindent
{\bf Proof of Theorem \ref{Lorentz}}. \emph{Part (i)}. By the monotonicity of the function in \eqref{mono1},
the couple of inequalities \eqref{red1} and \eqref{red2} in Theorem \ref{reduction}, applied with the norms appearing in inequality \eqref{mainlor}, is equivalent to 
 the couple of inequalities 
\begin{equation}\label{lor1}
\left\| s^{\frac{n-p}{\alpha p}-\frac {n-1}{\alpha}-\frac 1q}   \int _0^{s^{\frac
n\alpha}}\varphi(r) dr    \right\|_{L^{q}(0, \mu (\Omega) )} \leq
C 	\big\|s^{\frac 1p - \frac 1q}\varphi(s)\big\|_{L^{q}(0, \mathcal L^n (\Omega) )},
\end{equation}
and 
\begin{equation}\label{lor2}
\left\| \int _{s^{\frac n\alpha}}^\infty r^{-\frac{n-1}n}
\varphi(r) dr  \right\|_{L^{\frac{\alpha p}{n-p}, q}(0, \mu (\Omega) )} \leq
C \|\varphi\|_{L^{p,q}(0, \mathcal L^n (\Omega) )}
\end{equation}
for some constant $C$ and 
for every non-increasing function $\varphi: (0, \infty) \to [0,
\infty)$ with support in $(0, \mathcal L^n(\Omega))$. Moreover, by the monotonicity of the function in \eqref{mono2}
for any non-increasing function $\varphi$, inequality \eqref{red5} is equivalent 
to 
\begin{equation}\label{lor3}
\left\| s^{\frac{n-p}{\alpha p}-\frac {n-1}{\alpha}-\frac 1q}   \int _0^{s^{\frac
{n-1}\alpha}}\varphi(r) dr    \right\|_{L^{q}(0, \mu (\Omega) )} \leq
C \big\|s^{\frac {n-p}{p(n-1)} - \frac 1q}\varphi(s)\big\|_{L^{q}(0, \mathcal H^{n-1} (\partial \Omega) )},
\end{equation}
for some constant $C$ and  for every non-increasing function $\varphi: (0, \infty) \to [0,
\infty)$ with support in $(0, \hh(\partial \Omega))$. Inequalities \eqref{lor1} and \eqref{lor3} can be established via the criteria for weighted Hardy type inequalities in Lebesgue spaces  \cite[Section
1.3.2]{Mabook}. Inequality \eqref{lor2}  follows from the  proof    of \cite[Theorem 5.1]{CP_TAMS}.
With inequalities \eqref{lor1}--\eqref{lor3} in place, inequality \eqref{mainlor} is a consequence of Theorem \ref{reduction} and Remark \ref{remarkred}.
\\ \emph{Part (ii)}. The proof of inequality \eqref{BWlor} follows along the same  that of inequality \eqref{BWp}. The details are omitted, for brevity.
\\ \emph{Part (iii)}. Inequality \eqref{inflor} can be deduced from Theorem \ref{reduction} and Remark \ref{remarkred}, via weighted Hardy type inequalities in Lebesgue spaces  \cite[Section
1.3.2]{Mabook}.
\qed

\bigskip
\par
Our last result, contained in Theorem \ref{log}, consists of a set of inequalities for Orlicz norms  of power-logarithmic type, i.e. Zugmund norms, of the symmetric gradient. They extend the results of Theorems \ref{mainmeasure}, \ref{trudmeasure} and \ref{inf} in direction different from that of Theorem \ref{Lorentz}. It will be clear from the (sketched) proof that inequalities for more general Orlicz norms could be established via the same approach. Let us emphasize that the target norms are the strongest possible among all Orlicz spaces. Actually, they agree with the optimal Orlicz target norms appearing in parallel inequalities involving the full gradient. These inequalities are special cases of a result for arbitrary Orlicz-Sobolev spaces, established in \cite{Cianchisharp, Cianchibound} in a classical setting, namely for the Lebesgue measure and for functions vanishing on the boundary, or defined on  domains with some degree of regularity. The case of upper Ahlfors regular measures is considered in \cite{CPS}.
Let us point out that  an improvement of the conclusions of Theorem \ref{log} is however still possible, if more general rearrangement-invariant target norms are allowed. The proof, like that of  Theorem \ref{log}, relies upon Hardy type inequalities, with optimal rearrangement-invariant target spaces, from \cite{orlicztrace}.  Orlicz-Sobolev inequalities, with optimal rearrangement-invariant  norms, for the full gradient can be found in \cite{cianchiibero, cianchihigher} for the Lebesgue measure, and in \cite{CPS} for upper Ahlfors regular measures.  A discussion of these generalizations in the present framework is omitted, for brevity.

\begin{theorem}\label{log} {\bf [Zygmund--Sobolev inequalities]}
Let $\Omega$ be any  open set in $\rn$, $n \geq 2$. Let 
$\mu$ be a measure in $\Omega$ fulfilling \eqref{measure} for some exponent
$\alpha \in (n-1, n]$ and for some constant $C_\mu$.  
Assume that $\mathcal L^n (\Omega)< \infty$, $\mu (\Omega)<
\infty$ and $ \hh (\partial \Omega) < \infty$. 
\\ (i) Assume that $1 < p <
n$ and $\sigma \in \mathbb R$.    Then
 there
 exists a constant $C=C(n, p, \sigma, \alpha, C_\mu, \mathcal L^n (\Omega), \mu(\Omega), \hh (\partial \Omega))$ such that
\begin{align}\label{mainlog}
\|\bu\|_{
L^{\frac{p\alpha}{n-p}}(\log L)^{\frac{\sigma \alpha}{n-p}}
(\Omega, \mu)} & \leq
C \Big(\|\mathcal E \bu\|_{L^{p}(\log L)^\sigma(\Omega )} +   
\|\bu\|_{L^{\frac{p(n-1)}{n-p}}(\log L)^{\frac{\sigma (n-1)}{n-p}}(\partial \Omega )}\Big)
\end{align}
for every $u \in E^{1}L^{p}(\log L)^\sigma (\Omega) \cap C_b(\overline \Omega)$.
\\ (ii) Assume that $p=n$ and $\sigma < n-1$. 
Then there
 exists a constant  $C=C(n,  \sigma, \alpha, C_\mu, \mathcal L^n (\Omega), \mu(\Omega), \hh (\partial \Omega))$ such that
 \begin{align}\label{explog}
\|\bu\|_{\exp L^{\frac{n}{n-1-\sigma}}(\Omega, \mu)}  \leq
C\Big( \|\mathcal E \bu\|_{L^{n}(\log L)^\sigma(\Omega )}
  +
 \|\bu\|_{\exp L^{\frac{n}{n-1-\sigma}}(\partial \Omega)}\Big)
\end{align}
for every $u \in E^{1}L^{n}(\log L)^\sigma (\Omega) \cap C_b(\overline \Omega)$.
\\ (iii) Assume that $p=n$ and $\sigma = n-1$. 
Then there
 exists a constant  $C=C(n, \alpha, C_\mu, \mathcal L^n (\Omega), \mu(\Omega), \hh (\partial \Omega))$ such that
 \begin{align}\label{expexplog}
\|\bu\|_{\exp \exp L^{\frac{n}{n-1}}(\Omega, \mu)}  \leq
C\Big( \|\mathcal E \bu\|_{L^{n}(\log L)^{n-1}(\Omega )}
  +
 \|\bu\|_{\exp \exp L^{\frac{n}{n-1}}(\partial \Omega)}\Big)
\end{align}
for every $u \in E^{1}L^{n}(\log L)^{n-1}(\Omega) \cap C_b(\overline \Omega)$.
\\ (iv) Assume that either $p=n$ and $\sigma >n-1$, or $p>n$ and $\sigma \in \mathbb R$.
Then there
 exists a constant $C=C(n, p, \sigma, C_\mu, \mathcal L^n (\Omega), \mu(\Omega), \hh (\partial \Omega))$ such that
\begin{align}\label{inflog}
\|\bu\|_{L^\infty (\Omega )} & \leq C \Big(\|\mathcal E \bu\|_{L^{p}(\log L)^\sigma(\Omega )}
 +     \|\bu\|_{L^\infty(\partial \Omega)}\Big)
\end{align}
$u \in E^{1}L^{p}(\log L)^\sigma (\Omega) \cap C_b(\overline \Omega)$.
\end{theorem}
\medskip
\par\noindent
{\bf Proof, sketched}. By Theorem \ref{reduction} and Remark \ref{remarkred}, the proof is reduced to showing  the validity of the inequalities:
\begin{equation}\label{log1}
\left\| s^{-\frac {n-1}{\alpha}}\int _0^{s^{\frac
n\alpha}}\varphi(r) dr \right\|_{L^B(0, \infty )} \leq C
\|\varphi\|_{L^A(0, \infty )},
\end{equation}
\begin{equation}\label{log2}
\left\| \int _{s^{\frac n\alpha}}^\infty r^{-\frac{n-1}n}
\varphi(r) dr \right\|_{L^B(0, \infty )} \leq C
\|\varphi\|_{L^A(0, \infty )},
\end{equation}
for some constant $C$ and 
for every non-increasing function $\varphi: (0, \infty) \to [0,
\infty)$ with support in $(0, \mathcal L^n(\Omega))$, and
\begin{equation}\label{log3}
\left\| s^{-\frac {n-1}{\alpha}}\int _0^{s^{\frac {n-1}\alpha}}
\varphi(r)\, dr \right\|_{L^B(0, \infty )} \leq C
\|\varphi\|_{L^D(0, \infty )},
\end{equation}
for some constant $C$ and   for every non-increasing function $\varphi: (0, \infty) \to [0,
\infty)$ with support in $(0, \hh(\partial \Omega))$. Here, $A$, $B$, and $D$ are Young functions which yield the proper Orlicz spaces appearing in inequalities \eqref{mainlog}--\eqref{inflog}. Notice that, owing to the assumption that $\mathcal L^n (\Omega)< \infty$, $\mu (\Omega)<
\infty$ and $ \hh (\partial \Omega) < \infty$, only the behavior near infinity of the functions $A$, $B$, and $D$ is relevant here. Inequality \eqref{log2} is a special case of inequality (3.14) of \cite[Theorem 3.5]{orlicztrace}. The same theorem can be exploited to deal with inequalities \eqref{log1} and \eqref{log3}. Indeed, an H\"older type inequality in Orlicz spaces and Fubini's theorem ensure that
\begin{align}\label{log4}
& \sup _{\varphi \in L^A(0, \infty )} \frac {\left\|  s^{-\frac{n-1}\alpha} \int _0^{s^{\frac n\alpha}} 
\varphi(r) dr \right\|_{L^B(0, \infty )}}{\|\varphi\|_{L^A(0, \infty )}}
\approx
 \sup _{\varphi \in L^A(0, \infty )}\sup _{\psi \in L^{\widetilde B}(0, \infty)}  \frac {\int_0^\infty \psi(s) s^{-\frac{n-1}\alpha} \int _0^{s^{\frac n\alpha}} 
\varphi(r) dr\, ds}{\|\psi \|_{L^{\widetilde B}(0, \infty)}\|\varphi\|_{L^A(0, \infty )}}
\\
\nonumber & \approx \sup _{\psi \in L^{\widetilde B}(0, \infty)}  \sup _{\varphi \in L^A(0, \infty )} \frac {\int_0^\infty  
\varphi(r)  \int _ {r^{\frac \alpha n}}^\infty \psi(s)s^{-\frac{n-1}\alpha}  \, ds\,dr }{\|\psi \|_{L^{\widetilde B}(0, \infty)}\|\varphi\|_{L^A(0, \infty )}}
\approx \sup _{\psi \in L^{\widetilde B}(0, \infty)}\frac {\bigg\|  \int  _{r^{\frac \alpha n}}^\infty \psi(s)s^{-\frac{n-1}\alpha}  \, ds\bigg\|_{L^{\widetilde A}(0, \infty )}}{\|\psi \|_{L^{\widetilde B}(0, \infty)}}\,.
\end{align}
Here, $\widetilde A$ and $\widetilde B$ denote the Young conjugates of $A$ and $B$, and the relation $\lq\lq \approx"$ between two expressions means that they are bounded by each other, up to absolute multiplicative constants. Recall that $\widetilde A(t) = \sup\{st - A(s): s\geq 0\}$ for $t \geq 0$.
Therefore, inequality \eqref{log1} is equivalent to 
\begin{align}\label{log5}
\bigg\|  \int  _{r^{\frac \alpha n}}^\infty \psi(s)s^{-\frac{n-1}\alpha}  \, ds\bigg\|_{L^{\widetilde A}(0, \infty )}\leq \|\psi \|_{L^{\widetilde B}(0, \infty)}\,.
\end{align}
Similarly, inequality \eqref{log3} is equivalent to 
\begin{align}\label{log6}
\bigg\|  \int  _{r^{\frac \alpha {n-1}}}^\infty \psi(s)s^{-\frac{n-1}\alpha}  \, ds\bigg\|_{L^{\widetilde D}(0, \infty )}\leq \|\psi \|_{L^{\widetilde B}(0, \infty)}\,.
\end{align}
Both inequalities \eqref{log5} and \eqref{log6} can be established via \cite[Theorem 3.5]{orlicztrace}.
\qed

\section{Compliance with Ethical Standards}\label{conflicts}

\smallskip
\par\noindent 
{\bf Funding}. This research was partly funded by:   
\\ (i) Research Project of the
Italian Ministry of University and Research (MIUR) Prin 2015 \lq\lq  Partial differential equations and related analytic-geometric inequalities"  (grant number 2015HY8JCC);    
\\ (ii) GNAMPA   of the Italian INdAM - National Institute of High Mathematics (grant number not available);   
\\  (iii)   RUDN University Program 5-100.

\smallskip
\par\noindent
{\bf Conflict of Interest}. The authors declare that they have no conflict of interest.

%

%
%
%

\end{document}